\documentclass{amsart}[11pt]
\usepackage{mathrsfs, amsmath,amssymb}
\newtheorem{teo}{Theorem}
\newtheorem{pro}{Proposition}
\newtheorem{lem}{Lemma}
\newtheorem{cor}{Corollary}
\newtheorem*{rem}{Remark}

\title{On the spectral distribution of the free Jacobi process}

\author[N. Demni]{Nizar Demni}
\address{IRMAR, Universit\'e de Rennes 1\\ Campus de
Beaulieu\\ 35042 Rennes cedex\\ France}
\email{nizar.demni@univ-rennes1.fr}

\author[T. Hamdi]{Tarek Hamdi}
\address{IPEST, Universit\'e de Carthage\\ 2078 La Marsa \\ Tunisie}
\email{tarek.hamdi@ipest.rnu.tn}

\author[T. Hmidi]{Taoufik Hmidi}
\address{IRMAR, Universit\'e de Rennes 1\\ Campus de
Beaulieu\\ 35042 Rennes cedex\\ France}
\email{thmidi@univ-rennes1.fr}

\begin{document}
\maketitle

\begin{abstract}
In this paper, we are interested in the free Jacobi process starting at the unit of the compressed probability space where it takes values and associated with the parameter values $\lambda=1, \theta =1/2$. Firstly, we derive a time-dependent recurrence equation for the moments of the process (valid for any starting point and all parameter values). Secondly, we transform this equation to a nonlinear partial differential one for the moment generating function that we solve when $\lambda = 1, \theta =1/2$. The obtained solution together with tricky computations lead to an explicit expression of the moments which shows that the free Jacobi process  is distributed at any time $t$ as $(1/4)(2+Y_{2t}+Y_{2t}^{\star})$ where $Y$ is a free unitary Brownian motion. This expression is recovered relying on enumeration techniques after proving that if $a$ is a symmetric Bernoulli random variable which is free from $\{Y, Y^{\star}\}$, then the distributions of $Y_{2t}$ and that of $aY_taY_t^{\star}$ coincide. We close the exposition by investigating the spectral distribution associated with general sets of parameter values: $\lambda=1, \theta \in (0,1)$ and $\lambda \in (0,2), \theta = 1/2$. 
\end{abstract}

\section{Reminder and motivation}
The one dimensional Jacobi process is a two parameters-dependent stationary random motion valued in the interval $[0,1]$ and its equilibrium measure is given by the Beta distribution (\cite{Dem}). For special values of the parameters, this process can be realized as norms of projections of a spherical Brownian motion in a finite-dimensional Euclidean space onto spheres of lower dimensions (\cite{Bakry}). When extended to matrix spaces, this realization gives rise to matrix Jacobi processes as radial parts of corners of Brownian motions in the orthogonal and the unitary groups, and the corresponding equilibrium measure is the matrix-variate Beta distribution (\cite{Collins}, \cite{Dou}). Motivated by free probability theory, the free Jacobi process was realized in \cite{Demni} as the rescaled limit (in the sense of noncommutative moments) of the complex matrix Jacobi process. 
Its equilibrium distribution was determined in \cite{CapCas}, \cite{Collins} and \cite{Demni}. By the virtue of asymptotic freeness of independent unitarily-invariant random matrices and constant ones whose distributions converge, the matrix model of the free Jacobi process suggests the following abstract definition (see \cite{Demni} for details). Consider a non commutative W$^{\star}$-probability space $(\mathscr{A}, \tau)$, that is a von Neumann algebra $\mathscr{A}$ with unit ${\bf 1}$ and endowed with a faithful normalized trace $\tau$. Let $\theta \in (0,1), \lambda > 0$ such that $0 < \lambda \theta < 1$ and take two projections $P, Q$ in $\mathscr{A}$ such that 
\begin{itemize}
\item $\tau(P) = \lambda \theta,  \tau(Q) = \theta$, 
\item $PQ = QP = P$ ($P \leq Q$) if $\lambda \leq 1$ and $PQ=QP = Q$ ($P \geq Q$) otherwise. 
\end{itemize}
Let also $Y$ be a free unitary Brownian motion in $\mathscr{A}$ (\cite{Biane}) and assume
\begin{equation*}
\{P,Q\}, \,\, \{Y, Y^{\star}\} 
\end{equation*}
are free families in $\mathscr{A}$. Then 
\begin{equation*}
J_t := PY_tQY_t^{\star}P
\end{equation*}
defines the free Jacobi process of parameters $(\lambda, \theta)$ and one easily sees from 
\begin{equation*}
P-J_t = PY_t({\bf 1} - Q)Y_t^{\star}P 
\end{equation*}
that $P-J$ is a free Jacobi process too of parameters $(\lambda\theta/(1-\theta), 1-\theta)$. However, the stochastic analysis of $J$ performed in \cite{Demni} requires that $J_0$ and $P-J_0$ are injective operators in the compressed probability space 
\begin{equation*}
(P\mathscr{A}P, \frac{1}{\tau(P)} \tau).   
\end{equation*}
For that reason, the free Jacobi process studied in \cite{Demni} and denoted $J$ as well is driven by a free unitary Brownian motion starting at $Z$, where $Z$ is a unitary operator in $\mathscr{A}$ such that 
\begin{equation*}
\{P,Q\}, \,\{Z, Z^{\star}\}\,, \{Y, Y^{\star}\} 
\end{equation*}
are free families in $\mathscr{A}$. Actually, the operator $Z$ has to ensure that      
\begin{equation*}
T := \inf\{t > 0, \, J_t\, \textrm{or} \, P-J_t \, \, \textrm{is non injective}\}
\end{equation*}
is positive and the issue of our stochastic analysis is the free stochastic differential equation 
\begin{equation}\label{FSDE}
dJ_t = \sqrt{\lambda \theta} \sqrt{J_t}dW_t\sqrt{P-J_t} + \sqrt{\lambda \theta}\sqrt{J_t}dW_t^{\star} \sqrt{P-J_t} + (\theta P-J_t)dt, \quad 0 \leq t < T,
\end{equation}
where $W$ is a $P\mathscr{A}P$-complex free Brownian motion. This equation was the key ingredient for investigating spectral properties of $J_t$. Namely, let $m_n(t): = \tau(J_t^n)/\tau(P)$ be the moments of $J_t$ in $P\mathscr{A}P$ then for any $n \geq 1$ and up to time $T$ (\cite{Demni})
\begin{align}\label{Mom}
\partial_t m_n(t) & =  - nm_n(t)  + \theta nm_{n-1}(t)  + \lambda \theta n\sum_{k=0}^{n-2}m_{n-k-1}(t)(m_k(t) - m_{k+1}(t)) 
\end{align}
where the sum in the RHS is taken to be empty when $n=1$ and $m_0(t) =1, m_n(0) = \tau((PZQZ^{\star})^n)/\tau(P)$. Equation \eqref{Mom} may be thought of as a weak version of \eqref{FSDE} and since $J_t$ exists at any time $t$ and for any starting point $J_0$ then it is very likely that \eqref{Mom} remains valid for any initial data $m_n(0), n\geq 1$ (that is independently from $Z$). In this paper, we shall confirm the validity of \eqref{Mom} relying on the time-dependent recurrence equation derived in \cite{BenL}, Theorem 3.4 for 
\begin{equation*}
\tau(a_1Y_ta_2Y_t^{\star}\dots a_{2n-1}Y_ta_{2n}Y_t^{\star}),
\end{equation*}
where $\{a_1, \dots, a_{2n}\}$ is a $2n$-tuple of random variables forming a $\star$-free family with $Y$. More precisely, we first specialize this theorem to $a_k = P$ when $k$ is odd and $a_k = Q$ when $k$ is even, then infer that the computations keep unchanged when we rather consider $a_k = Z^{\star}PZ$ when $k$ is odd and $Z$ is any unitary operator in $\mathscr{A}$ which is $\star$-free from $Y$ and from $\{P,Q\}$. Indeed, $\tau$ is a trace, $Z^{\star}PZ$ and $P$ have the same moments and both families $\{Z^{\star} PZ, Q\}$ and $\{Y, Y^{\star}\}$ are free. Doing so gives us the benefit of studying the free Jacobi process $PYQY^{\star}P$ which was discarded in \cite{Demni} and which is related to Voiculescu's liberation process associated with free projections (\cite{Voiculescu}). However, observe that for any $\lambda$, the sequence $v_n(t) := \lambda m_n(t), n \geq 1$ satisfies 
\begin{align*}
\partial_t v_n(t) & =  - nv_n(t)  + \theta n v_{n-1}(t)  + \theta n\sum_{k=0}^{n-2}v_{n-k-1}(t)(v_k(t) - v_{k+1}(t)) 
\end{align*}
with initial data $v_n(0) = \lambda m_n(0), v_0(t) = \lambda m_0(t)$. Thus, we shall primarily consider the value $\lambda=1$ and focus on the corresponding free Jacobi process. However, due to the high nonlinearity of \eqref{Mom}, we further restrict our attention to the value $\theta =1/2$ and prove that the spectral distribution of $J_t$ in the compressed space, say $\mu_t$, fits the distribution of the random variable 
\begin{equation*}
\frac{1}{4}[Y_{2t} + Y_{2t}^{\star} + 2{\bf 1}]
\end{equation*}
in $(\mathscr{A}, \tau)$. In particular $\mu_t$ is absolutely continuous with respect to Lebesgue measure on the real line and its support fill in the interval $(0,1)$ at time $t=2$. Actually, our description follows from a closed formula for $m_n(t), n \geq 1$, namely: 
\begin{align}\label{Form}
m_n(t) =  \frac{1}{2^{2n}}  \binom{2n}{n} + \frac{1}{2^{2n-1}} \sum_{k=1}^n  \binom{2n}{n-k}\frac{1}{k} L_{k-1}^{1}(2kt) e^{-kt}
\end{align}
where $L_n^{1}$ is the $n$-th Laguerre polynomial of index $1$ (\cite{Rainville}). Formula \eqref{Form} is in turn obtained after solving a nonlinear partial differential equation (hereafter p.d.e) for the moment generating function of $\mu_t$:  
\begin{equation*}
M_t(z) = \sum_{n \geq 0} m_n(t)z^n  \quad |z| < 1.
\end{equation*}
In fact, by the virtue of \eqref{Mom}, this p.d.e admits a unique solution in the class of analytic functions around zero. Nevertheless, the binomial coefficients present in \eqref{Form} suggests that this formula might be derived based on enumeration techniques. To this aim, we start by noting that  when $\lambda =1$, the faithfulness of $\tau$ and $P \leq Q$ entail $P=Q$: indeed $\tau(P) = \tau(Q)$ so that 
\begin{equation*}
\tau((P-Q)^2) = 2\tau(P) - 2\tau(PQ) = 0. 
\end{equation*} 
Next, if further $\theta = 1/2$ then $P = (1+a)/2$ where $a \in \mathscr{A}$ is a self-adjoint symmetric Bernoulli random variable and we shall derive the following expansion: 
\begin{equation*}
2\tau(({\bf 1}+a)Y_t({\bf 1}+a)Y_t^{\star})^n) =  \binom{2n}{n} + \sum_{k=1}^n \binom{2n}{n-k}\tau((aY_taY_t^{\star})^k).
\end{equation*}
By the way, we are grateful to A. Nica and to C. Sattlecker for their help in completing the proof of this derivation. Finally, we recover \eqref{Form} after proving that $Y_{2t}$ and $v_t := aY_taY_t^{\star}$ have the same distribution. The enumeration techniques not only provides another elegant proof of the description of $\mu_t$, but opens the way to investigate the spectral distribution of $J_t$ for arbitrary $\theta \in (0,1)$ and $\lambda = 1$. More precisely, the above expansion of $\tau(({\bf 1}+a)Y_t({\bf 1}+a)Y_t^{\star})^n)$ comes in with the additional term $2^{2n-1}(2\theta-1)$. However one has to describe the spectral distribution of $aY_taY_t^{\star}$ which does not seem to be easy. As to the parameter values $\lambda \in (0,2), \theta = 1/2$, we already noticed above that they correspond to a deformation of the initial data $m_n(0), m_0(t)$. However we should consider the parameter values $\lambda \in (0,1]$ and $\lambda \in [1,2)$ for fixed $\theta = 1/2$ separately, since $PQ = P \Rightarrow m_n(0) = 1$ and $PQ = Q \Rightarrow m_n(0) = 1/\lambda$ respectively. Nonetheless, note that the moments $(m_n(t))_n$ corresponding to $\lambda \in [1,2)$ can be recovered from the ones corresponding to $\lambda \in (0,1]$. Indeed, if $P$ is a projection of continuous rank $\tau(P) = \lambda /2$ then ${\bf 1} - P$ is a projection of continuous rank $(2-\lambda)/2$. Moreover, since $Y_tQY_t^{\star}$ is idempotent and since $\tau$ is tracial then we easily derive  
\begin{equation*}
\tau[(({\bf 1} - P)Y_tQY_t^{\star})^n] = \tau(Q) + \sum_{k=1}^n (-1)^{k} \binom{n}{k} \tau[(PY_tQY_t^{\star})^k]. 
\end{equation*}
Accordingly, we only need to consider the parameter values $\lambda \in (0,1], \theta = 1/2$ and we shall prove that the moment generating function of the spectral distribution of $J_t$ splits into two parts. The first `principal' part is the sum of  the moment generating function of the stationary distribution (i.e. the spectral distribution in the compressed space of $PUQU^{\star}P$ where $U$ is a Haar unitary variable in $\mathscr{A}$ which is $\star$-free from $\{P,Q\}$) and of the generating function of a linear combination of the moments of $Y_{2\lambda t/(2-\lambda)}$. The second part is an analytic function $u_t$ around the origin whose first and second derivatives vanish there. Unfortunately, we do not succeed to provide a description of the spectral distribution of the free Jacobi process through known random variables as we already did for $\lambda = 1, \theta = 1/2$.  

\section{A recurrence time-dependent equation}
In this section, we prove that
\begin{pro}
The moments of the free Jacobi process 
\begin{equation*}
m_n(t):= \frac{1}{\tau(P)}\tau(PZY_tQY_t^{\star}Z^{\star}P)^n), \, n \geq 1
\end{equation*}
satisfy \eqref{Mom} for any unitary $Z$ that is $\star$-free from $Y$ and from $\{P,Q\}$.  
\end{pro}
{\it Proof}: As explained in the introductory part, computations are given in the special case $Z = {\bf 1}$. Moreover, we shall equivalently prove that the sequence 
\begin{equation*}
r_n(t):= \tau(J_t^n) = \tau(P) m_n(t), n \geq 1
\end{equation*}
satisfies   
\begin{align}\label{RR}
\partial_t r_n(t) & =  - nr_n(t)  + n\theta r_{n-1}(t)  + n\sum_{k=0}^{n-2}r_{n-k-1}(t)(r_k(t) - r_{k+1}(t))
\end{align}
with $r_0(t) = \tau(P) = \lambda \theta$. To proceed, we recall Theorem 3.4 in \cite{BenL}: 

\begin{teo}\label{Theo}
Let $n \geq 1$ and define
\begin{equation*}
f_{2n}(a_1, \dots, a_{2n}, t) := e^{nt} \tau(a_1Y_ta_2Y_t^{\star}\dots a_{2n-1}Y_ta_{2n}Y_t^{\star})
\end{equation*}
where $\{a_1, \dots, a_{2n}\} \in \mathscr{A}$ is $\star$-free with $Y$. Set $f_0(A,t) := \tau(A)$ for any $A \in \mathscr{A}$ then 
\begin{align*}
&\partial_t  f_{2n}(a_1, \dots, a_{2n}, t) = - \sum_{\substack{1 \leq k < l \leq 2n \\ l - k \equiv 0[2]}}f_{2n - (l-k)}(a_1, \dots, a_k, a_{l+1}, \dots, a_{2n}, t)f_{l-k}(a_{k+1}, \dots, a_l,t) +
\\&  e^t\sum_{\substack{1 \leq k < l \leq 2n \\ l-k -1\equiv 0[2]}}f_{2n - (l-k)-1}(a_1, \dots, a_{k-1}, a_k a_{l+1}, a_{l+2}, \dots, a_{2n}, t)f_{l-k-1}(a_la_{k+1}, a_{k+2}, \dots, a_{l-1},t)
\end{align*}
\end{teo}

Now, we specialize Theorem \ref{Theo} to $a_k = P$ if $k$ is odd and $a_k = Q$ otherwise. Then the result is straightforward when $n=1$ and is even stated in \cite{BenL}, p.923. So let $n \geq 2$ and note that both indices $k,l$ in the first (respectively second) sum in Theorem \ref{Theo} have the same (respectively different) parity, therefore $k$ and $l+1$ in the second sum have the same parity and so do $l$ and $k+1$. Accordingly, the first sum does not contain terms $f_0(\cdot, t)$ while the second does: they correspond to indices $l=2n, k=1$ and to $l=k+1, 1 \leq k \leq 2n-1$. Since $P$ and $Q$ are idempotent and since $\tau$ is a trace, then the contribution of indices $k=1, l=2n$ is $\tau(P)f_{2n-2}(P,Q, \dots, P, Q)$ while that of $l=k+1, 1 \leq k \leq 2n-1$ is
\begin{equation*}
[n\tau(Q) + (n-1)\tau(P)]f_{2n-2}(P, Q, \dots, P, Q, t) 
\end{equation*}
 respectively. Thus, both contributions sum up to 
 \begin{equation}\label{C0}
 n(\tau(Q)+\tau(P))f_{2n-2}(P, Q, \dots, P, Q, t).
 \end{equation}
Next we write $l = k+2s+1$ for integer positive values of $s$ and distinguish $n = 2$ and $n \geq 3$. If $n = 2$ then there is no additional term in the second sum, while if $n \geq 3$ we separate $k=1$ and $2 \leq  k  \leq 2n-3$. By the same properties of $P, Q, \tau$ mentioned above, the contribution of indices $k=1, l=2s+2$ is 
 \begin{equation}\label{C1}
\sum_{s=1}^{n-2}f_{2(n- s-1)}(P,Q, \dots, P, Q, t)f_{2s}(P, Q, \dots, P,Q, t). 
 \end{equation}
 For the remaining values of $2 \leq k \leq 2n-3$, we distinguish even and odd ones: the contribution of indices $k=2j, 1 \leq j \leq n-2, l = 2j+2s+1$ is 
 \begin{equation*}
 \sum_{j=1}^{n-2}\sum_{s=1}^{n-j-1} f_{2(n- s-1)}(P,Q, \dots, P, Q, t)f_{2s}(P, Q, \dots, P,Q, t) 
 \end{equation*}
while that of $k=2j+1, 1 \leq j \leq n-2, l=2s+2j+2$ is also 
 \begin{equation*}
 \sum_{j=1}^{n-2}\sum_{s=1}^{n-j-1} f_{2(n- s-1)}(P,Q, \dots, P, Q, t)f_{2s}(P, Q, \dots, P,Q, t).  
 \end{equation*}
By rearranging the terms in both obtained double sums, we see that the contribution of indices $2 \leq k \leq 2n-3, l = 2k+2s+1$ is 
\begin{equation*}
2 \sum_{s=1}^{n-2}(n-s-1) f_{2(n- s-1)}(P,Q, \dots, P, Q, t)f_{2s}(P, Q, \dots, P,Q, t)
\end{equation*}
which simplify after the index change $s \mapsto n-s-1$ to 
\begin{equation}\label{C2}
(n-1)\sum_{s=1}^{n-2} f_{2(n- s-1)}(P,Q, \dots, P, Q, t)f_{2s}(P, Q, \dots, P,Q, t). 
\end{equation}
 
Similarly we consider the first sum in Theorem \ref{Theo} and write $l=k+2s, 1 \leq k \leq 2n-2, 1 \leq s \leq n-[(k+1)/2], n \geq 2$. Then it contributes to
\begin{equation}\label{C3}
n \sum_{s=1}^{n-1} f_{2(n- s)}(P,Q, \dots, P, Q, t)f_{2s}(P, Q, \dots, P,Q, t). 
\end{equation}
Now we recall that $f_{2n}(P, Q, \dots, P, Q, t) = e^{nt}r_n(t)$ and take into account the exponential factor in front of the second sum of Theorem \ref{Theo}. If $n \geq 3$ then \eqref{C2} and \eqref{C3} sum up to
\begin{equation*}
e^{nt}[-nr_{n-1}(t)r_1(t) - \sum_{s=1}^{n-2}r_{n-s-1}(t)r_{s}(t) + n\sum_{s=1}^{n-2}r_{s}(t)[r_{n-s-1}(t) - r_{n-s}(t)].
 \end{equation*}
With regard to \eqref{C0}, \eqref{C1}, the whole contribution of the RHS of Theorem \ref{Theo} is  
\begin{equation*}
e^{nt}\{n\tau(Q) + n[\tau(P) - r_0(t)] r_{n-1}(t) + n\sum_{s=1}^{n-1}r_{s}(t)[r_{n-s-1}(t) - r_{n-s}(t)]\}.
\end{equation*}
But $\partial_tf_{2n}(t) = e^{nt}[\partial_t r_nt + nr_n(t)]$ together with $r_0(t)= \tau(P)$ show that \eqref{RR} holds for any $t > 0$ and any $n \geq 3$. If $n=2$ then the whole contribution is given by \eqref{C0} and \eqref{C3} leading to \eqref{RR} as well.  $\hfill \blacksquare$

\section{The case $\theta = 1/2, \lambda = 1$}
This section is devoted to the description of the spectral distribution $\mu_t$ of $J_t$ when $\lambda = 1, \theta = 1/2$ and $J_0 = P$. A major step towards it is the following result:  
 \begin{pro}
Let $L_k^1$ be the $k$-th Laguerre polynomial and let
\begin{equation*}
\rho_t(z) := \sum_{k=1}^{\infty} \frac{1}{k}L_{k-1}^1(kt)z^k, \, |z| < 1,
\end{equation*}
be the unique solution of (\cite{Biane}, \cite{Rains}) 
\begin{equation*}
\partial_t\rho_t + \frac{z}{2}\partial_z\rho_t^2 = 0, \quad \rho_0(z) = \frac{z}{1-z}
\end{equation*} 
in the class of analytic function around zero. Then the moment generating function $M_t$ of $\mu_t$ is given by
 \begin{equation*}
 M_t(z) = \frac{1}{\sqrt{1-z}}\left\{1+ 2 \rho_{2t} \left(\frac{ze^{-t}}{(\sqrt{1-z}+1)^{2}}\right)\right\}, \, |z| < 1.
 \end{equation*}
\end{pro}
 
{\it Proof}: Before coming through computations, we point out that 
$$
\alpha: \, z \mapsto \frac{z}{(1+\sqrt{1-z})^2}
$$
maps the open unit disc into itself. Indeed, the following expansion holds (\cite{Rainville} p.70) 
\begin{equation*}
\frac{1}{(1+\sqrt{1-z})^2} = \frac{1}{4} \sum_{k \geq 0} \frac{(1)_n(3/2)_n}{(3)_n n!} z^n, \, |z| < 1
\end{equation*} 
and is still convergent for $z=1$ by Gauss Theorem (\cite{Rainville}, p.49). Hence the expression of $M_t(z)$ given in the proposition makes sense for all $|z| < 1$. Now, let $|z| > 1$ then similar computations leading to Proposition 7.1. in \cite{Demni} shows that 
\begin{equation*}
G_t(z) :=  \frac{1}{z}\sum_{n=0}^{\infty}\frac{m_n(t)}{z^n},
\end{equation*}
satisfies the p.d.e. for all $(\lambda, \theta)$
\begin{equation*}
\partial_tG_t = \partial_z\left\{[(1-2\lambda \theta)z - \theta(1-\lambda)]G_t + \lambda \theta z(z-1)G_t^2 \right\}
\end{equation*}
with initial value $G_0(z) = 1/(z-1)$. This p.d.e. simplifies when one substitutes $\lambda = 1, \theta = 1/2$ to  
\begin{equation*}
\partial_tG_t = \frac{1}{2} \partial_z\left\{z(z-1)G_t^2\right\},\, G_0(z) = \frac{1}{z-1}
\end{equation*}  
or equivalently for $|z| < 1$ 
\begin{equation*}
\partial_tM_t = -\frac{z}{2} \partial_z\left\{(1-z)M_t^2\right\}, \, M_0(z) = \frac{1}{1-z}.
\end{equation*}  
 Set $S_t(z):=\sqrt{1-z}M_t(z)-1,$ then
 $$
 \partial_t S_t+z\sqrt{1-z}\partial_zS_t+\frac12z\sqrt{1-z}\partial_zS_t^2=0
 $$
 with the initial data 
 \begin{equation*}
 S_0(z) = \frac{1}{\sqrt{1-z}} - 1.
 \end{equation*}
Next note that $\alpha$ is invertible in a neighborhood of zero with inverse function given by 
\begin{equation*}
\alpha^{-1}(z) = \frac{4z}{(1+z)^2}.
\end{equation*}
Keeping in mind that $|\alpha(z)| \leq |z| < 1$ in the open unit disc, then $\alpha$ extends to a biholomorphic map from the open unit disc onto its image. Moreover,
$$
z\sqrt{1-z}\alpha^\prime(z)=\alpha(z).
$$
Hence $F_t: z \mapsto S_t(\alpha^{-1}(z))$ and $w_t: z \mapsto  F_t(e^t z)$ satisfy
$$
\partial_tF_t+z\partial_zF_t+\frac12 z\partial_z F_t^2=0
$$
and
$$
\partial_tw_t+\frac12 z\partial_z w_t^2=0.
$$
Finally one easily checks from 
\begin{equation*}
1+ \alpha(z) = \frac{2}{1+\sqrt{1-z}},\quad 1-\alpha(z) = \frac{2\sqrt{1-z}}{1+\sqrt{1-z}}
\end{equation*}
that
\begin{equation*}
w_0(\alpha(z)) = S_0(z) = \frac{2\alpha(z)}{1-\alpha(z)} = 2\rho_0(\alpha(z))\, \Leftrightarrow \, w_0(z) = 2\rho_0(z).
\end{equation*}
Since $2\rho_{2t}$ and $w_t$ satisfy the same partial differential equation as $\rho_t$, the proposition is proved. $\hfill \blacksquare$

%\begin{rem}
%It is interesting to check whether the function $\alpha$ yields a subordination property for the free unitary Brownian motion (\cite{Biane0}, Theorem 3.5). 
%\end{rem}

The moments $m_n(t)$ are given by 
\begin{cor}
%\begin{equation*}
%Q_n(t) := \frac{(-1)^n}{2^{2n+1}}\frac{1}{n+1}L_n^{1}(4(n+1)t), n \geq 0.
%\end{equation*}
For any $n \geq 1$ and any $t \geq 0$
\begin{align*}
m_n(t) & %\frac{(1/2)_n}{n!}  + \sum_{k=1}^n \frac{(-1)^{k-1}}{2^{2(n-k)}} \binom{2n}{n-k}Q_{k-1}(t/2)e^{-kt}
 =  \frac{1}{2^{2n}}\binom{2n}{n}  + \frac{1}{2^{2n-1}} \sum_{k=1}^n  \binom{2n}{n-k}\frac{1}{k} L_{k-1}^1(2kt) e^{-kt}.
\end{align*}
\end{cor}
%\begin{equation*}
%m_n(\infty) = \frac{\Gamma(n+1/2)}{\sqrt{\pi}n!} = \
%\end{equation*}

{\it Proof}: First of all the generalized binomial Theorem yields
\begin{equation*}
\frac{1}{\sqrt{1-z}} = \sum_{n = 0}^{\infty}\frac{(1/2)_n}{n!}z^n, \, |z| < 1
\end{equation*}
where $(1/2)_n = \Gamma(n+1/2)/\Gamma(1/2)$ is the Pochhammer symbol. But Legendre duplication formula (\cite{Erd}) shows that 
\begin{equation*}
\frac{(1/2)_n}{n!} = \frac{1}{2^{2n}}\binom{2n}{n}.
\end{equation*}
 
Now let ${}_2F_1$ be the Gauss hypergeometric function then for any $k \geq 1$ (\cite{Rainville} p.70)
\begin{align*}
\frac{2^{2k}}{\sqrt{1-z}}\frac{1}{(\sqrt{1-z}+1)^{2k}} &= {}_2F_1\left(k+\frac{1}{2}, k+1, 2k+1, z\right) 
\\& = \sum_{n = 0}^{\infty}\frac{(k+1/2)_n(k+1)_n}{(2k+1)_n} \frac{z^n}{n!}
\end{align*}
%\frac{2^{2k} \Gamma(k+1/2)\Gamma(k+1)}{\sqrt{\pi} \Gamma(2k+1)} 
Using Legendre duplication formula again, we derive 
\begin{equation*}
(2n+2k)! = \frac{2^{2n+2k}}{\sqrt{\pi}} \Gamma(n+k+1/2)\Gamma(n+k+1) 
\end{equation*}
which yields 
\begin{align*}
\sum_{n =0}^{\infty} \binom{2n+2k}{n} \frac{z^n}{2^{2n}} &= \frac{2^{2k} \Gamma(k+1/2)\Gamma(k+1)}{\sqrt{\pi} \Gamma(2k+1)} \sum_{n = 0}^{\infty}\frac{(k+1/2)_n(k+1)_n}{(2k+1)_n} \frac{z^n}{n!}
\\& = \frac{2^{2k}}{\sqrt{1-z}}\frac{1}{(\sqrt{1-z}+1)^{2k}}.
\end{align*}
As a result 
\begin{align*}
\frac{1}{\sqrt{1-z}}\rho_{2t} \left(\frac{ze^{-t}}{(\sqrt{1-z}+1)^{2}}\right) &= \frac{1}{\sqrt{1-z}}\sum_{k=1}^{\infty} \frac{1}{k} L_{k-1}^1(2kt)\left[\frac{ze^{-t}}{(\sqrt{1-z}+1)^{2}}\right]^k
\\& = \sum_{k=1}^{\infty} \frac{1}{k} L_{k-1}^1(2kt) e^{-kt} \sum_{n =0}^{\infty} \binom{2n+2k}{n} \frac{z^{n+k}}{2^{2n+2k}}
\\& = \sum_{k=1}^{\infty} \frac{1}{k} L_{k-1}^1(2kt) e^{-kt} \sum_{n =k}^{\infty} \binom{2n}{n-k} \frac{z^{n}}{2^{2n}}. 
\end{align*}
The Corollary is proved. $\hfill \blacksquare$

We are now ready to give the description of $\mu_t$: 
\begin{cor}
When $\theta=1/2, \lambda = 1$, the free Jacobi process starting at $P$ is distributed in the compressed probability space $(P\mathscr{A}P, 2\tau)$ as the random variable
\begin{equation*}
\frac{Y_{2t}^{-1}}{4}(1+Y_{2t})^2  = \frac{1}{4}[Y_{2t}^{-1} + 2 + Y_{2t}]
 \end{equation*}
in $\mathscr{A}$. Consequently, $\mu_t$ is absolutely continuous with respect to Lebesgue measure on the real line and its density is given by
\begin{equation*}
2 \frac{k_{2t}(e^{2i\arccos(\sqrt{x})})}{\sqrt{x(1-x)}} {\bf 1}_{[0,1]}(x) dx. 
\end{equation*}
where $k_t$ is the density of the spectral distribution of the free unitary Brownian motion. In particular, the support of $\mu_t$ is the whole interval $[0,1]$ for any time $t \geq 2$. 
\end{cor}

{\it Proof}: Let  $h_n(t) = \tau(Y_t^n), n \geq 1$ be the moments of $Y_t$ in $(\mathscr{A}, \tau)$, then (\cite{Biane})
\begin{equation*}
h_n(t) = \frac{e^{-nt/2}}{n} L_{n-1}^1(nt),
\end{equation*}
But since $Y$ and $Y^{\star} = Y^{-1}$ have the same distribution then $h_n = h_{-n}, n \in \mathbb{Z}$. As a result
\begin{align*}
m_n(t) &= \frac{\displaystyle \binom{2n}{n}}{2^{2n}}  + \frac{2}{2^{2n}} \sum_{k=1}^n  \binom{2n}{n-k}h_{k}(2t)
 \\& =  \frac{1}{2^{2n}} \left\{\binom{2n}{n}  + \sum_{k=1}^n \binom{2n}{n-k}h_{k}(2t) + \sum_{k=1}^{n} \binom{2n}{n+k}h_{k}(2t)\right\}
\\& =  \frac{1}{2^{2n}} \left\{\binom{2n}{n}  + \sum_{k=1}^n \binom{2n}{n-k}h_{k}(2t) + \sum_{k=-n}^{-1} \binom{2n}{n-k}h_{-k}(2t)\right\}
\\& = \frac{1}{2^{2n}}  \sum_{k=-n}^n \binom{2n}{n-k}h_{k}(2t) 
\\& = \frac{1}{2^{2n}}  \sum_{k=0}^{2n} \binom{2n}{2n-k}\tau(Y_{2t}^{k-n})
\\& =  \frac{1}{2^{2n}}\tau(Y_{2t}^{-n}(1+Y_{2t})^{2n}).
 \end{align*}
Finally, the support of $\mu_t$ is entirely determined by the one of $k_{2t}$ (see \cite{Biane1}).

\begin{rem}
When $\lambda = 1, \theta = 1/2$  \eqref{Mom} takes the form 
\begin{align*}
\partial_tm_n(t) + \frac{n}{2}m_n(t) = \frac{n}{2} \sum_{k=0}^{n-1}m_{n-1-k}(t)[m_k(t) - m_{k+1}(t)].
\end{align*}
Thus since the moments of $Y_t$ converge as $t \infty$ to those of a Haar unitary random variable, then the moments 
\begin{equation*}
m_n(\infty) = \frac{1}{2^{2n}} \binom{2n}{n} := m_n
\end{equation*}
satisfy
\begin{equation}\label{Comb}
m_n = \sum_{k=0}^{n-1}m_{n-k-1} [m_{k} - m_{k+1}]. 
\end{equation}
But if $C_k:= (1/(k+1))\displaystyle \binom{2k}{k}$ is the $k$-th Catalan number (\cite{NS}) then Legendre duplication formula entails 
\begin{equation*}
m_k - m_{k+1} = \frac{1}{2\sqrt {\pi}} \frac{\Gamma(k+1/2)}{(k+1)!} = \frac{1}{2^{2k+1}} C_k. 
\end{equation*}
As a result
\begin{align*}
\sum_{k=0}^{n-1}m_{n-k-1} [m_{k} - m_{k+1}] = \frac{2}{2^{2n}} \sum_{k=0}^{n-1}(n-k)C_{n-k-1}C_k = \frac{n+1}{2^{2n}} \sum_{k=0}^{n-1} C_{n-k-1}C_k. 
\end{align*}
Consequently, \eqref{Comb} is nothing else but the recurrence relation for Catalan numbers: 
\begin{equation*}
C_n = \sum_{k=0}^{n-1} C_{n-1-k}C_k. 
\end{equation*} 
\end{rem}

\section{Enumerative derivation of the moments} 
Recall that if $\lambda = 1$ then the faithfulness of $\tau$ forces $P=Q$. If further $\theta = 1/2$ then $P := (1+a)/2$ where $a=a^{\star} \in \mathscr{A}$ is distributed according to 
\begin{equation*}
\frac{1}{2}(\delta_1+\delta_{-1}).
\end{equation*}
In the sequel, we shall derive \eqref{Form} relying on enumeration techniques. More precisely, 
\begin{pro}
Let $a$ be a self-adjoint random variable in $\mathscr{A}$ distributed according to 
\begin{equation*}
\frac{1}{2}(\delta_1+\delta_{-1})
\end{equation*}
and assume $a$ and $Y$ are $\star$-free. Then for any $n \geq 1$
\begin{equation*}
\tau[((1+a)Y_t(1+a)Y_t^{\star})^n] =  \frac{1}{2}\binom{2n}{n}  + \sum_{k=1}^n  \binom{2n}{n-k}\tau((aY_taY_t^{\star})^k).
\end{equation*}
\end{pro}

{\it Proof}: let $b := Y_taY_t^{\star}$ then 
\begin{equation*}
(1+a)Y_t(1+a)Y_t^{\star} = (1+a)(1+b)
\end{equation*}
therefore $[(1+a)(1+b)]^n$ consists of words formed by letters picked in the alphabet $\{a,b\}$ subject to the cancellations $\{a^{2k}=b^{2k}={\bf 1}, a^{2k+1}= a, b^{2k+1}= b\}$. Moreover, we claim that those formed by an odd number of letters have zero expectation. Indeed, any such word has a reduced (after taking into account of cancellations and that there is no relation between $a$ and $b$) expression either $aba\cdots aba$ or $bab\cdots bab$. The claim then follows from the trace property of $\tau$ and from $\tau(a) = \tau(b) = 0$. As a result, we need to enumerate for each $n \geq 1$ words $\{(ab)^k, 1 \leq k \leq n, (ba)^k, 1 \leq k \leq n-1\}$ (see below for the constant term $k=0$). To this end, let $c(n,k), d(n,k), e(n,k)$ be the number of words $(ab)^k, (ab)^ka, (ba)^k$ respectively. Then the expansion 
\begin{equation*}
[(1+a)(1+b)]^n = [(1+a)(1+b)]^{n-1}(1+a+b+ab)
\end{equation*}
and the observation 
\begin{equation*}
(ab)^k = [(ab)^k] {\bf 1} = [(ab)^ka] a = [(ab)^{k-1}a]b = [(ab)^{k-1}] (ab)
\end{equation*}
give the recurrence relation 
\begin{equation*}
c(n,k) = c(n-1,k) + d(n-1,k) + d(n-1,k-1) + c(n-1, k-1), \, n \geq 2, k \geq 1
\end{equation*}  
and 
\begin{equation*}
c(1,k) = \delta_{k0} + \delta_{k1} 
\end{equation*}
where $\delta_{ij}$ is the Kronecker symbol. In a similar fashion, 
\begin{equation*}
(ab)^ka = [(ab)^ka] {\bf 1} = [(ab)^k]a = [(ab)^{k+1}]b = [(ab)^{k+1}a] (ab)
\end{equation*}
gives the recurrence relation
\begin{equation*}
d(n,k) = d(n-1,k) + c(n-1,k) + c(n-1,k+1) + d(n-1,k+1), \, n\geq 2, k \geq 0,
\end{equation*}
with 
\begin{equation*}
d(1,k) = \delta_{k0},
\end{equation*}
while 
\begin{equation*}
[(1+a)(1+b)]^n = (1+a+b+ab)[(1+a)(1+b)]^{n-1}
\end{equation*}
together with 
\begin{equation*}
(ba)^k =  {\bf 1}(ba)^k = a[(ab)^ka] = b[(ab)^{k-1}a] = ab[(ba)^{k+1}]
\end{equation*}
give the recurrence relation 
\begin{equation*}
e(n,k) = e(n-1,k) + d(n-1,k) + d(n-1,k-1) + e(n-1,k+1), \, n\geq 3, k \geq 1,
\end{equation*}
with 
\begin{equation*}
e(2,k) = 3\delta_{k0} + \delta_{k1}, \, e(1,k) = \delta_{k0}.
\end{equation*}
Now, we compute $c(n,0)$ which correspond to the constant term (not depending on time $t$) in the expansion 
\begin{equation*}
\tau[((1+a)(1+b))^n].  
\end{equation*}
To proceed, we argue that up to a factor, $c(n,0)$ is the $n$-th moment $m_n(\infty)$ of the stationary distribution of the free Jacobi process 
\begin{equation*}
2\frac{c(n,0)}{2^{2n}} = \frac{1}{2^{2n}}\binom{2n}{n}\, \Rightarrow c(n,0) = \frac{1}{2}\binom{2n}{n}.
\end{equation*}
Indeed, for fixed $n \geq 1$
\begin{equation*}
\lim_{t \rightarrow \infty} \tau[(aY_taY_t^{\star})^n] = \tau[(aUaU^{\star})^n]
\end{equation*}
where $U$ is a Haar unitary distributed random variable. But $a$ and $UaU^{\star}$ are free in $\mathscr{A}$ and since $\tau(a) = \tau(UaU^{\star}) = 0$ then the very definition of freeness implies that $\tau[(aUaU^{\star})^n] = 0$ for any $n \geq 1$. We can also compute $c(n,0)$ by considering the von Neumann algebra generated by $\{a,b\}$ endowed with the state that assigns the value $1$ to the unit and vanishes otherwise. More precisely, this algebra is the free product of the von Neumann algebras generated by $\{a\}$ and by $\{b\}$ since the latters may be realized as von Neumann algebras of the cyclic group of order two so that Theorem 1.6.3. in \cite{DNV} applies. Therefore $a$ and $b$ are free there and it is obvious that they are symmetric Bernoulli random variable with respect to the given state. It follows that $c(n,0)$ is the constant term in the $n$-th moment of the two-fold free multiplicative convolution 
\begin{equation*}
\frac{1}{2}(\delta_0 + \delta_2) \boxtimes \frac{1}{2}(\delta_0 + \delta_2)
\end{equation*}
which may be computed using the $S$-transform (\cite{NS}). More precisely, the $S$-transform of this free multiplicative convolution reads
\begin{equation*}
\frac{(z+1)^2}{(1+2z)^2}
\end{equation*}
near the origin. Consequently, the corresponding inverse moment generating is given by 
\begin{equation*}
\frac{z(z+1)}{(1+2z)^2} = \frac{1}{4}\left[1-\frac{1}{(1+2z)^2}\right]
\end{equation*}
whence we deduce the inverse moment generating of the free multiplicative convolution
\begin{equation*}
\frac{1}{2}\left[\frac{1}{\sqrt{1-4z}} - 1\right] = \frac{1}{2} \sum_{n =1}^{\infty}\binom{2n}n z^n.
\end{equation*}

Based on the initial values $c(n,0), d(1,0)$, we proceed by mutual induction (on $k$ for fixed $n$ then on $n$) and use the elementary identity 
\begin{equation*}
\binom{n}{k} + \binom{n}{k-1} = \binom{n+1}{k}
\end{equation*}
in order to prove that 
\begin{eqnarray*}
c(n,k) &=& \binom{2n-1}{n-k}\\
d(n,k) &=& \binom{2n-1}{n-k-1} = e(n,k).
\end{eqnarray*}
Note by passing that $d(n, k-1) = c(n, k)$ which agrees with the recurrence relations for $d(n,k)$ and $e(n,k)$. Finally 
\begin{equation*}
c(n,k)+e(n,k) = \binom{2n-1}{n-k} + \binom{2n-1}{n-k-1} = \binom{2n}{n-k}. \hfill \blacksquare  
\end{equation*}

In order to recover \eqref{Form}, it suffices to observe that a
\begin{lem}
The unitary random variable $aY_taY_t^{\star}$ is distributed as $Y_{2t}$.
\end{lem}
{\it Proof}: since the distribution of a unitary random variable is entirely determined by its moments, we shall prove that the sequence defined by
\begin{equation*}
s_n(t):= e^{nt} \tau((aY_taY_t^{\star})^n), \, n\geq 1, 
\end{equation*}
is the (unique) solution of 
\begin{equation*}
\partial_ts_n(t) = -n\sum_{k=1}^{n-1}s_{n-k}(t) s_k(t), \quad n \geq 2, \, s_1(t) = 1,
\end{equation*}
 then infer from \cite{Rains} (Theorem 4 p. 669) that 
\begin{equation*}
s_n(t) = \frac{1}{n}L_{n-1}^1(2nt).
\end{equation*}
To proceed, we specialize Theorem \ref{Theo} to $a_k = a$ for all $k$ and notice that the second sum vanishes since $a^2 = {\bf 1}$ and since $\tau(a) = 0$. The contribution of the first sum is easily seen to be
\begin{equation*}
-n\sum_{k=1}^{n-1}s_{n-k}(t) s_k(t)
\end{equation*}
and the value of $s_1(t)=1$ is readily derived from \cite{BenL} p. 923 using $\tau(a) = 0, \tau(a^2) = 1$. The lemma is proved. $\hfill \blacksquare$ \\

\section{Further developments: general parameter values}
In this section, we investigate the spectral distribution of $J_t$ for more general sets of parameter values: $\lambda =1, \theta \in (0,1)$, $\lambda \in (0,2), \theta = 1/2$. 

\subsection{Parameter values $\lambda = 1, \theta \in (0,1)$} 
As one easily realizes, the enumeration techniques used to expand 
\begin{equation*}
\tau[((1+a)(1+b))^n]
\end{equation*}
remain valid when $\lambda =1, \theta \in (0,1)$. Indeed, $P= (1+a)/2$ where the spectral distribution of $a$ is 
\begin{equation*}
\theta \delta_1 + (1-\theta)\delta_{-1},
\end{equation*}
and we need to take into account the contribution of words with odd number of letters. By the trace property of $\tau$ and the relations $a^2=b^2 = {\bf 1}$, this contribution is $\tau(a) = (2\theta-1)$ up to a positive integer say $c(n)$ therefore 
\begin{equation*}
\tau[((1+a)Y_t(1+a)Y_t^{\star})^n] = \frac{1}{2}\binom{2n}{n} + \sum_{k=1}^n \binom{2n}{n-k}\tau((aY_taY_t^{\star})^k) + (2\theta - 1)c(n).  
\end{equation*}
The coefficient $c(n)$ is easily computed when $\theta \in (0,1)\setminus \{1/2\}$ by letting $t=0$ and using $(1+a)^2 = 4P, a^2 = {\bf 1}$: 
\begin{equation*}
4^n \theta = 2^{2n-1} + (2\theta - 1) c(n) \, \Rightarrow \, c(n) = 2^{2n-1}.
\end{equation*}
However, it is not easy to determine the distribution of $aY_taY_t^{\star}$ for general $\theta \in (0,1)$: indeed let again $s_n(t) = e^{nt} \tau((aY_taY_t^{\star})^n), \, n\geq 1$ then 
\begin{eqnarray*}
\partial_t s_n(t) & =&  -n\sum_{k=1}^{n-1}s_{n-k}(t) s_k(t) + e^{nt}[2n(2\theta-1) + (n-1)(n-2)(2\theta-1)^2], \, n \geq 2, \\
s_1(t) & = & 1 + (2\theta-1)^2(e^t-1).
\end{eqnarray*}

\subsection{Parameter values $\lambda \in (0,1], \theta = 1/2$}
In this paragraph, we shall prove the following decomposition: 
\begin{pro}
Let 
\begin{equation*}
M_t(z):= \sum_{n = 0}^{\infty} m_n(t)z^n, \quad t \geq 0, |z| < 1, 
\end{equation*}
be the moment generating function of the free Jacobi process associated with the parameter values $\lambda \in (0,1), \theta =1/2$ (we omit the dependence on $\lambda$ for sake of clarity). 
Let also $M_{\infty}$ be the moment generating function of the stationary distribution of the free Jacobi process associated with the same parameter values (the distribution of $PUQU^{\star}P$ when $U \in \mathscr{A}$ is Haar distributed and is $\star$-free from $\{P,Q\}$). Then 
\begin{equation*}
M_t(z) = M_{\infty}(z) + \frac{2}{2-\lambda} \rho_{2\lambda t/(2-\lambda)} [(2-\lambda)e^{-t}\alpha(z)] + u_t(z)
\end{equation*}
where $u_t$ is analytic around the origin whose first and second derivatives vanish, and
\begin{equation*}
\lim_{\lambda \rightarrow 1}u_t(z) = 0.
\end{equation*}
\end{pro}

{\it Proof}: Recall that for general parameter values $(\lambda, \theta)$, the function 
\begin{equation*}
(t,z) \mapsto G_t(z) = \frac{1}{z}\sum_{n \geq 0} \frac{m_n(t)}{z^n} = \frac{1}{z}M_t\left(\frac{1}{z}\right)
\end{equation*} 
defined for $t > 0$ and $|z| > 1$, is a solution of the p.d.e. 
\begin{equation*}
\partial_tG_t = \partial_z\left\{[(1-2\lambda \theta)z - \theta(1-\lambda)]G_t + \lambda \theta z(z-1)G_t^2 \right\}
\end{equation*}
with initial value $G_0(z) = 1/(z-1)$. Recall also from \cite{Demni} that 
\begin{equation*}
G_{\infty}(z) := \frac{(2-r)z + (1/\lambda - 1) + \sqrt{r^2z^2 - Bz + C^2}}{2z(z-1)}
\end{equation*}
where $r = (1/\lambda \theta), \, B = 2(r + (r-2)/\lambda), C = (1-1/\lambda)$, is the Cauchy-Stieltjes transform of the stationary distribution of the free Jacobi process. Then 
\begin{equation*}
\partial_z\left\{[(1-2\lambda \theta)z - \theta(1-\lambda)]G_{\infty} + \lambda \theta z(z-1)G_{\infty}^2 \right\} = 0
\end{equation*}
as can be readily checked from the expression of $G_{\infty}$. In particular, for $\theta = 1/2, \lambda \in (0,1]$ the p.d.e   
\begin{equation*}
\partial_tG_t = \frac{1}{2}\partial_z\left\{(1-\lambda)(2z - 1)G_t + \lambda z(z-1)G_t^2 \right\}
\end{equation*}
holds, while
\begin{equation*}
G_{\infty}(z) := \frac{(\lambda-1)(2z-1) + \sqrt{4z^2 - 4z + (1-\lambda)^2}}{2\lambda z(z-1)} = \frac{1}{z}M_{\infty}\left(\frac{1}{z}\right)
\end{equation*}
satisfies 
\begin{equation*}
\partial_z\left\{(1-\lambda)(2z - 1)G_{\infty} + \lambda z(z-1)G_{\infty}^2 \right\} = 0.
\end{equation*}
Set 
\begin{equation*}
H_t := G_t - G_{\infty}
\end{equation*}
then
\begin{equation*}
\partial_t H_t = \frac{1}{2} \partial_z \left\{\lambda z(z-1)H_t^2 + \sqrt{4z^2 - 4z + (1-\lambda)^2}H_t\right\}.
\end{equation*}
In particular, if $z \in (0,1]$ and  
\begin{equation*}
S_t(z) := \frac{1}{z}H_t \left(\frac{1}{z}\right) = M_t(z) - M_{\infty}(z)
\end{equation*}
then 
\begin{equation*}
\partial_tS_t = -\frac{z}{2}\partial_z \left\{\lambda (1-z)S_t^2 + \sqrt{4 - 4z + (1-\lambda)^2z^2}S_t\right\}. 
\end{equation*}
Now denote
\begin{equation*}
v_t(z) := \frac{2}{2-\lambda} \rho_{2\lambda t/(2-\lambda)} [(2-\lambda)e^{-t}\alpha(z)]
\end{equation*}
where we recall that 
\begin{equation*}
\rho_t(z) = \sum_{k \geq 1}\frac{1}{k}L_{k-1}^1(kt)z^k, \quad \alpha(z) = \frac{z}{(1+\sqrt{1-z})^2}. 
\end{equation*}
Then the p.d.e 
\begin{equation*}
\partial_t \rho_t + \frac{z}{2}\partial_z\rho_t^2 = 0
\end{equation*}
together with
\begin{equation*}
\frac{\alpha'(z)}{\alpha(z)} = \frac{1}{z\sqrt{1-z}}
\end{equation*}
show that
\begin{equation*}
\partial_tv_t (z)+ \frac{\lambda}{2} z\sqrt{1-z}\partial_zv_t^2 (z) = -2e^{-t}\alpha(z)(\partial_z \rho_t)[(2-\lambda)e^{-t}\alpha(z)].
\end{equation*}
Accordingly, the function  
\begin{equation*}
(t,z) \mapsto u_t(z):= S_t(z) -\frac{1}{\sqrt{1-z}} v_t(z) = \frac{1}{z}H_t \left(\frac{1}{z}\right) - \frac{1}{\sqrt{1-z}} v_t(z)
\end{equation*}
satisfies 
\begin{align*}
\partial_tu_t &= -\frac{z}{2}\partial_z \left\{\lambda(1-z)u_t^2 + 2\lambda\sqrt{1-z}u_tv_t + \sqrt{4 - 4z + (1-\lambda)^2z^2}u_t + \right. 
\\& \left. \sqrt{4 + (1-\lambda)^2\frac{z^2}{1-z}}v_t\right\}(z) - \frac{1}{\sqrt{1-z}}\left[\partial_tv_t + \frac{\lambda}{2}z\sqrt{1-z}\partial_z v_t^2\right](z)
\\& = -\frac{z}{2}\partial_z \left\{\lambda(1-z)u_t^2 + 2\lambda\sqrt{1-z}u_tv_t + \sqrt{4 - 4z + (1-\lambda)^2z^2}u_t + \right. 
\\& \left. \sqrt{4 + (1-\lambda)^2\frac{z^2}{1-z}}v_t\right\}(z) + \frac{2e^{-t}\alpha(z)}{\sqrt{1-z}}(\partial_z \rho_t)[(2-\lambda)e^{-t}\alpha(z)]
\\& = -\frac{z}{2}\partial_z \left\{\lambda(1-z)u_t^2 + 2\lambda\sqrt{1-z}u_tv_t + \sqrt{4 - 4z + (1-\lambda)^2z^2}u_t + \right. 
\\& \left. \sqrt{4 + (1-\lambda)^2\frac{z^2}{1-z}}v_t\right\}(z) + 2ze^{-t}\alpha'(z)(\partial_z \rho_t)[(2-\lambda)e^{-t}\alpha(z)].
\end{align*}

Set
\begin{equation*}
r_t(z):= \sqrt{4 + (1-\lambda)^2\frac{z^2}{1-z}}
\end{equation*}
then 
\begin{align*}
-\frac{z}{2}\partial_z(r_tv_t)(z) = -\frac{(1-\lambda)^2}{4r_t(z)}\frac{z^2(2-z)}{(1-z)^2}v_t(z) - z r_t(z)\alpha'(z)e^{-t}(\partial_z \rho_t)[(2-\lambda)e^{-t}\alpha(z)].
\end{align*}
Rearranging terms, one gets 
\begin{align*}
\partial_tu_t &= -\frac{z}{2}\partial_z \left\{\lambda(1-z)u_t^2 + 2\sqrt{1-z}u_tv_t + \sqrt{4 - 4z + (1-\lambda)^2z^2}u_t  \right\} 
\\& +\frac{(1-\lambda)^2}{4r_t(z)}\frac{z^2(z-2)}{(1-z)^2}v_t(z) + ze^{-t}(r_t(z) - 2)\alpha'(z)(\partial_z \rho_t)[(2-\lambda)e^{-t}\alpha(z)].
\end{align*}
Since 
\begin{equation*}
\lim_{\lambda \rightarrow 1} r_t(z)= 2 
\end{equation*}
then one can see that 
\begin{equation*}
Z_t : =\frac{(1-\lambda)^2}{4r_t(z)}\frac{z^2(z-2)}{(1-z)^2}v_t(z) + ze^{-t}(r_t(z) - 2)\alpha'(z)(\partial_z \rho_t)[(2-\lambda)e^{-t}\alpha(z)] \rightarrow 0 
\end{equation*}
as $\lambda \rightarrow 1$.  As a matter of fact, one can expand
\begin{equation*}
Z_t = (1-\lambda)\sum_{n = 1}^{\infty}d_n(t)z^n
\end{equation*}
for some coefficients $d_n(t)$ depending on $t,\lambda$. Now write 
\begin{equation*}
u_t(z) = \sum_{n=1}^{\infty}c_n(t)z^n
\end{equation*}
and recall the expansion 
\begin{equation*}
\sqrt{1-z} = \sum_{n = 0}^{\infty}\frac{-(2n)!}{(2n-1)(2^{n}n!)^2}z^n := \sum_{n = 0}^{\infty}\beta_nz^n.
\end{equation*}
Then for $\lambda \neq 1$
\begin{equation*}
\sqrt{4 - 4z + (1-\lambda)^2z^2} = \sum_{n = 0}^{\infty}\gamma_nz^n 
\end{equation*}
where 
\begin{equation*}
\gamma_n = 2\sum_{k = 0}^n \frac{\beta_k\beta_{n-k}}{z_1^kz_2^k} z^k
\end{equation*}  
and
\begin{equation*}
z_1= \frac{2 + 2\sqrt{\lambda(2-\lambda)}}{(1-\lambda)^2}, \quad z_2= \frac{2 - 2\sqrt{\lambda(2-\lambda)}}{(1-\lambda)^2}.
\end{equation*}
Besides, 
\begin{eqnarray*}
\sqrt{4-4z+(1-\lambda)^2z^2}u_t(z) &=& \sum_{n = 1}^{\infty}\left( \sum_{k=1}^{n} c_k(t)\gamma_{n-k}\right)z^n \\
\frac{1}{\sqrt{1-z}}v_t(z)u_t(z) &=& \sum_{n = 2}^{\infty}\left( \sum_{k=1}^{n-1} c_k(t)\psi_{n-k}(t)\right)z^n 
\end{eqnarray*}
where we set 
\begin{equation*}
\frac{1}{\sqrt{1-z}}v_t(z) = \sum_{n=1}^{\infty}\psi_n(t)z^n
\end{equation*}
with
\begin{equation*}
\psi_n(t) := \frac{1}{2^{2n-1}} \sum_{k=1}^n\frac{1}{k}(2-\lambda)^{k-1}L_{k-1}^1\left(\frac{2\lambda t}{2-\lambda} k\right)e^{-kt}, \, n \geq 1.
\end{equation*}
Consequently, 
\begin{eqnarray*}
c'_1(t)&=& -c_1(t)+(1-\lambda)d_1(t).
\end{eqnarray*}

But 
\begin{align}
 u_0(z) &= \frac{1}{z}G_0\left(\frac{1}{z}\right) - \frac{1}{z}G_{\infty}\left(\frac{1}{z}\right) - \frac{1}{\sqrt{1-z}}v_0(z) \nonumber \\
 &=\frac{1}{1-z}-\left(\frac{(1-\lambda)(z-2)+\sqrt{4-4z+(1-\lambda)^2z^2}}{2\lambda(1-z)} \right) -\frac{1}{\sqrt{1-z}}v_0(z)\nonumber \\
&= \left( \frac{1}{\lambda}-\frac{1-\lambda}{2\lambda}z -\frac{1}{2\lambda}\sqrt{4-4z+(1-\lambda)^2z^2}\right)\frac{1}{1-z}-\frac{1}{\sqrt{1-z}}v_0(z) \nonumber \\
&=\left( \frac{1}{\lambda}-\frac{1-\lambda}{2\lambda}z -\frac{1}{2\lambda}\sum_{n\geq0}\gamma_n z^n\right)\sum_{n\geq0} z^n-\sum_{n\geq1}\psi_n(0) z^n \nonumber \\
&=\sum_{n=1}^{\infty}\left( \frac{1+\lambda}{2\lambda}-\frac{1}{2\lambda}\sum_{k=0}^n\gamma_k-\psi_n(0)\right) z^n \label{Coeff}.
\end{align}
In particular 
\begin{equation*}
c_1(0)=\frac{1+\lambda}{2\lambda}-\frac{1}{2\lambda}-\frac{1}{2}=0
\end{equation*}
which implies that 
\begin{equation*}
c_1(t) = (1-\lambda) \left(\int_0^td_1(s)e^s ds\right) e^{-t}.
\end{equation*}
Since $d_1(t) = 0$ then $c_1(t) = 0$ and 
\begin{equation*}
u_t(z) = \sum_{n=2}^{\infty} c_n(t)z^n. 
\end{equation*}
Accordingly 
\begin{equation*}
c'_2(t) = -2c_2(t)+(1-\lambda)d_2(t),
\end{equation*}
and since 
\begin{equation*}
c_2(0)=\frac{1+\lambda}{2\lambda}-\frac{1}{2\lambda}(1-\frac{\lambda(2-\lambda)}{4})-\frac{1}{8}(4+2-\lambda)=0
\end{equation*}
then 
\begin{equation*}
c_2(t) = (1-\lambda) \left(\int_0^td_2(s)e^{2s} ds\right) e^{-2t}.
\end{equation*}
But after some effort, one sees that $d_2(t) = 0$ yielding $c_2(t) = 0$ and 
\begin{equation*}
u_t(z) = \sum_{n=3}^{\infty} c_n(t)z^n. 
\end{equation*}
Similarly 
\begin{equation*}
c'_3(t) = -3c_3(t) + (1-\lambda)d_3(t) 
\end{equation*}
and 
\begin{align*}
c_3(0)%&=\frac{1+\lambda}{2\lambda}-\frac{1}{2\lambda}(1-\frac{\lambda(2-\lambda)}{4}-\frac{\lambda(2-\lambda)}{8})-\frac{15+6(2-\lambda)+(2-\lambda)^2}{32}\\
=-\frac{(1-\lambda)(3-\lambda)}{32}.
\end{align*}
It follows  that 
\begin{equation*}
c_3(t) = (1-\lambda) \left(-\frac{3-\lambda}{32}+\int_0^td_3(s)e^{3s}ds\right) e^{-3t}
\end{equation*}
and computer assisted computations show actually that 
\begin{equation*}
c_3(t) = - \frac{1-\lambda}{32}\left[2\lambda e^{-3t} + 3(1-\lambda)e^{-t}\right]. 
\end{equation*}

More generally
\begin{align}\label{General}
c'_n(t) &= - \frac{n}{2}\left( \sum_{k=3}^{n} c_k(t)\gamma_{n-k}  + 2\psi_1(t)c_{n-1}(t)+\sum_{k=3}^{n-2}\left[\lambda(c_{n-k}(t)-c_{n-1-k}(t)) \right. \right. \nonumber
\\& \left. \left.+2(\psi_{n-k}(t)-\psi_{n-1-k}(t))\right] c_k(t) \right) +(1-\lambda)d_n(t).
\end{align}
 for $n \geq 4$, where the last sum is taken to be zero when $n=4$. Finally, it remains to prove  that 
\begin{equation*}
\lim_{\lambda \rightarrow 1} c_n(t) = 0. 
\end{equation*}
To proceed, it is easy to note from \eqref{Coeff} that 
\begin{equation*}
c_n(0)=c_{n-1}(0) + k_n(\lambda)
\end{equation*}
for any  $n\geq 2$, where
\begin{eqnarray*}
 k_n(\lambda)&=&\psi_{n-1}(0)-\psi_n(0)-\frac{1}{2\lambda}\gamma_n\\
&=&-\frac{1}{2^{2n-1}}+\frac{1}{2^{2n-1}}\sum_{k =1}^{n-1}\left( 4\binom{2n-2}{n-1-k}-
\binom{2n}{n-k}\right)-\frac{1}{\lambda} \sum_{k=0}^n\frac{\beta_k\beta_{n-k}}{z_1^kz_2^{n-k}}.
\end{eqnarray*}
But
\begin{equation*}
0\leq \frac{1}{z_1^kz_2^{n-k}}=\frac{(1-\lambda)^{2k}(1+\sqrt{\lambda(2-\lambda)})^{n-2k}}{2^n} \leq 1
\end{equation*}
together
\begin{equation*}
\lim_{\lambda \rightarrow 1} \sqrt{4-4z + (1-\lambda)^2z^2}  = 2\sqrt{1-z}
\end{equation*}
entail
\begin{equation*}
\lim_{\lambda \rightarrow 1}\gamma_n=\lim_{\lambda \rightarrow 1} 2\sum_{k=0}^n\frac{\beta_k\beta_{n-k}}{z_1^kz_2^{n-k}}=\frac{-2(2n)!}{(2n-1)(2^{n}n!)^2}.
\end{equation*}
As a result
\begin{eqnarray*}
\lim_{\lambda \rightarrow 1} k_n(\lambda)%&=&
%-\frac{1}{2^{2n-1}}+\frac{1}{2^{2n-1}}\sum_{k =1}^{n-1}\left( 4\binom{2n-2}{n-1-k}-
%\binom{2n}{n-k}\right)-\lim_{\lambda \rightarrow 1}\frac{1}{\lambda} \sum_{k=0}^n\frac{\beta_k\beta_{n-k}}{\lambda_1^k\lambda_2^{n-k}}\\
&=&-\frac{1}{2^{2n-1}}+\frac{1}{2^{2n-1}}\sum_{k =1}^{n-1}\left( 4\binom{2n-2}{n-1-k}-
\binom{2n}{n-k}\right)+\frac{(2n)!}{(2n-1)(2^{n}n!)^2}\\
\end{eqnarray*}
and applying twice Pascal rule shows that
\begin{eqnarray*}
\binom{2n}{n-k}=\binom{2n-2}{n-k-2}+\binom{2n-2}{n-k}+2\binom{2n-2}{n-k-1}.
\end{eqnarray*}
Hence
\begin{equation*}
4\binom{2n-2}{n-1-k} - \binom{2n}{n-k}=\binom{2n-2}{n-k-1}-\binom{2n-2}{n-k-2}+\binom{2n-2}{n-k-1}-\binom{2n-2}{n-k}.
\end{equation*}
yielding
\begin{eqnarray*}
\lim_{\lambda \rightarrow 1}  k_n(\lambda) %&=&-\frac{1}{2^{2n-1}}+\frac{1}{2^{2n-1}}\left( \binom{2n-2}{n-2} - \binom{2n-2}{n-1}+1\right)+\frac{(2n)!}{(2n-1)(2^{n}n!)^2}\\
&=&\frac{1}{2^{2n-1}}\left( \binom{2n-2}{n-2}-\binom{2n-2}{n-1}\right)+\frac{(2n)!}{(2n-1)(2^{n}n!)^2} = 0. 
\end{eqnarray*}
As a matter fact, it suffices to use induction in order to prove that $c_n(0) \rightarrow 0, n \geq 4$ as $\lambda \rightarrow 1$. Once this is proved, it is clear from \eqref{General} that
\begin{equation*}
\lim_{\lambda \rightarrow 1} c_n(t)=0
\end{equation*}
for any $n \geq 4$ and any $t \geq 0$. $\hfill \blacksquare$

\subsection{Parameter values $\lambda \in (1,2), \theta =1/2$}
As mentioned in the introduction, the moments of the free Jacobi process associated with the parameter values $\lambda \in (1,2), \theta =1/2$ may be deduced from those corresponding to 
$\lambda \in (1,2), \theta =1/2$. Indeed
\begin{pro} 
\begin{equation*}
\tau[(({\bf 1}-P)Y_tQY_t^*)^n]=\tau(Q)+\sum_{k=1}^n(-1)^k\binom{n}{k} \tau((PY_tQY_t^*)^k).
\end{equation*}
\end{pro}

{\it Proof}: Let $n\geq1$ and define for $a_1, a_2, \dots, a_n \in \mathscr{A}$  
\begin{equation*}
h_n(a_1,\ldots,a_n,t):=\tau(a_1Y_tQY_t^*\ldots a_nY_tQY_t^*).
\end{equation*}
$h_n$ is obviously a linear map with respect to its $n$ first coordinates, thereby 
\begin{eqnarray*}
h_n({\bf 1}-P,\ldots, {\bf 1}-P,t) = \sum_{(a_{i_1},\ldots,a_{i_n})\in\{1,-P\}^n}h_n(a_{i_1},\cdots,a_{i_n},t)
\end{eqnarray*}
Since $Q$ is idempotent and $\tau$ is a trace, then one realizes that each $n$-tuples with $k, 1\leq k \leq n$ elements equal to $(-P)$ and $n-k$ elements equal to ${\bf 1}$ contributes to $(-1)^{k}\tau((PY_tQY_t^*)^k)$. They are even the only tuples giving this contribution since there is no relation between $P$ and $Y_tQY_t^{\star}$. Since their total contribution is 
\begin{equation*}
(-1)^k\binom{n}{k} \tau[(PY_tQY_t^*)^k],\quad 1\leq k \leq n.
\end{equation*}
while that of the $n$-tuple $({\bf 1},\dots, {\bf 1})$ is $\tau(Q)$, then the proposition is proved. $\hfill \blacksquare$

\end{document}